\documentclass[%draft
]{amsproc}
\usepackage{amsmath}
\usepackage{amssymb}
\usepackage{amsthm}
\usepackage{graphicx}

\RequirePackage[cp1251]{inputenc}
\RequirePackage[english,russian]{babel}
\RequirePackage{srcltx}

\theoremstyle{plain}

\theoremstyle{definition}

\numberwithin{equation}{section}

\begin{document}

\title[О подчинении линейных операторов]{О подчинении линейных операторов на лебеговых пространствах над $\mathbb R^d$\\
On subordination of linear operators on the Lebesgue\\ spaces over $\mathbb R^d$}

\author{Роальд М. Тригуб}


\date{}

\maketitle

\begin{abstract}
In this note a general approach is suggested for comparison of operators.
This is done by means of the Fourier transform of a measure. This approach is applied
to comparison of approximation properties of various summability methods of the Fourier integrals (I)
and to differential operators with constant coefficients (II).
\end{abstract}

\vskip4mm

Речь пойдёт о неравенствах вида
$$ \|\tilde{\Lambda}f\|_q\le c\|\Lambda f\|_p$$
и
$$ \|\tilde{\Lambda}f\|_q\le c(\|\Lambda_1 f\|_{p_1}+\|\Lambda_2 f\|_{p_2})$$
с некоторой константой $c$,не зависящей от функции $f$.

Здесь и далее $ \|f\|_p$ - $L_p$-норма, $p\ge 1,$ над евклидовым пространством $\mathbb R^d$ со
стандартным скалярным произведением. В качестве $L_\infty$ часто предполагают
пространство непрерывных ограниченных функций с $sup$-нормой модуля функции.

Основная идея состоит в попытке представить операторы, стоящие справа в неравенстве, в
виде свёртки с мерой и последующем применении неравенств для свёрток.

Через $\mu$ будем обозначать конечную на $\mathbb R^d$ борелевскую комплекснозначную
меру, а через $(\mathbb R^d,\Omega,\mu)$-пространство с мерой.

Для любого множества $E\in \Omega$ положим
$$var \mu (E)=|\mu|(E)=\sup_{E=\cup_1^\infty E_k}\sum |\mu(E_k)| \quad (E_k\in
\Omega,E_{k_1}\cap E_{k_2}=\O,k_1\ne k_2),
$$
$$
\widehat {d\mu}(y)=\int_{\mathbb R^d} e^{-i(x,y)}d\mu(x)
$$
и
$$
W=W(\mathbb R^d)=\{\widehat {d\mu},||\widehat {d\mu}||_W= |\mu| (\mathbb R^d)<\infty\}.
$$
Как следует из теоремы Радона-Никодима, существует функция $h_\mu$ такая, что
$|h_\mu(x)|=1$ почти всюду по мере $|\mu|$ и для любого $E\in \Omega$ и $f\in
L_1(E,\mu)$
$$
\int_E fd\mu=\int_E fh_\mu d|\mu|.
$$
См., напр., [1],гл.11.

 Если мера $\mu$ абсолютно непрерывна относительно меры Лебега на
$\mathbb R^d$ ($d\mu(x)=g(x)dx$, $g\in L_1(\mathbb R^d)$), то получаем идеал винеровской
банаховой алгебры $W$:
$$
W_0=W_0(\mathbb R^d)=\{f(x)=\hat g(x)=\int_{\mathbb
R^d}g(y)e^{-i(x,y)}dy,||f||_{W_0}=||g||_1<\infty\}.
$$
Отличие в поведении функций из $W$ и $W_0$ лишь около $\infty$. В частности, если $f\in
W$, $\lim_{|x|\to \infty} f(x)=0$ и вне некоторого куба функция $f$ ограниченной вариации
по Витали, то $f\in W_0 [2]$.

Обозначим через $W_1$ подалгебру $W$, которая получается из $W_0$ присоединением единицы.

См. обзорную статью $[3]$, содержащую необходимые и достаточные условия принадлежности
этим алгебрам.

Если теперь $\tilde {\Lambda}f=\Lambda f*d\mu$ - свёртка, т.е.
$$
\tilde {\Lambda}f(x)=\int_{\mathbb R^d}\Lambda f(x-y)d\mu(y),
$$
то для $p\in [1,\infty]$
$$
\|\tilde \Lambda f\|_p\le |\mu|(\mathbb R^d)\|\Lambda f\|_p.
$$
 В частности, при $d\mu(x)=g(x)dx$
 $$
\|\tilde \Lambda f\|_p\le \|g\|_1\|\Lambda f\|_p
$$
и при $g(x)\ge 0$ константу $\|g\|_1$ нельзя уменьшить, вообще говоря ($[4]$,гл.I,(4.2)).
А если при этом $g\in L_s$, $s>1,$ то в силу неравенства Юнга для свёрток (см., напр.,
[4],гл.I,(4.3)) $\Lambda f \in L_q,\frac 1q=\frac 1p +\frac 1s -1.$

Обозначим через $Z(\psi)$ множество нулей функции $\psi\in C(\mathbb R^d)$.

\underline{Лемма 1} (принцип сравнения)

Если $Z(\widehat {d\mu}_2)\subset Z(\widehat {d\mu}_1)$ и $\widehat {d\mu}_1=\widehat
{d\mu}_{1,2}\cdot \widehat {d\mu}_2$, то для любого $p\ge 1$
$$
\|f*d\mu_1\|_p\le K\|f*d\mu_2\|_p,
$$
где
$$
K=\inf_{\frac 00}\|\frac {\widehat {d\mu}_1}{\widehat {d\mu}_2}\|_W=\inf
|\mu_{1,2}|(\mathbb R^d)
$$
(нижняя грань относится к выбору значений $\widehat {d\mu}_{1,2}$ на $Z(\widehat
{d\mu}_2)$).

Доказательство.

Если $f\in L_1(\mathbb R^d)$, то и $f*d\mu \in L_1(\mathbb R^d)$, а
$$
\widehat {f*d\mu_1}=\hat f\cdot \widehat {d\mu}_1=\widehat {d\mu}_{1,2} \cdot (\hat
f\cdot \widehat {d\mu}_2)=\widehat {d\mu}_{1,2}\cdot \widehat {f*d\mu_2}.
$$
Если равны преобразования Фурье двух функций  из $L_1$, то функции совпадают почти всюду
(по мере Лебега). Поэтому
$$
f*d\mu_1=d\mu_{1,2}*(f*d\mu_2),\|f*d\mu_1\|_1\le |\mu_{1,2}|(\mathbb R^d)\|f*d\mu_2\|_1.
$$
Если $f\in L_p,p>1,$ то применяем доказанное равенство к функции $f_n$, совпадающей с $f$
при $|x|\le n$ и равной нулю при $|x|>n$:
$$
f_n*d\mu_1=d\mu_{1,2}*(f_n*d\mu_2,\|f_n*d\mu_1\|p\le |\mu_{1,2}|(\mathbb
R^d)\|(f_n*d\mu_2)\|_p.
$$
Осталось при $p<\infty$ в этом неравенстве перейти к пределу при $n\to
\infty$, учитывая, что
$$
\|(f-f_n)*d\mu\|_p\le \|f-f_n\|_p|\mu|\to 0,
$$
а при при $p=\infty$ в указанном равенстве перейти к пределу, используя теорему Лебега о
мажорируемой сходимости.

Лемма доказана. Эта лемма есть, по сути, в [5].

Свёртки - это мультипликаторы Фурье [4]. Приведенный принцип сравнения уточняется для
периодических функций и применяется к методам суммирования кратных рядов Фурье функций
на торе $\mathbb T^d$ ([2],[6],8.3).

\emph{I. Сравнение разных методов суммирования интегралов Фурье (на примере средних типа
Гаусса-Вейерштрасса).}

Пусть $\phi_\alpha(x)=e^{-|x|^\alpha},\alpha>0.$ Для функций $f\in L_1(\mathbb R^d)$ в
силу формулы умножения (см., напр., [4],гл.1)
$$
\frac {1}{(2\pi)^d}\int_{\mathbb R^d}\phi_\alpha(\epsilon y)\hat f(y)e^{i(x,y)}dy= \frac
{1}{(2\pi)^d}\int_{\mathbb R^d}\hat\phi_\alpha (y)f(x-\epsilon y)dy.
$$
При $\alpha=1$ и $\alpha=2$ сходимость этих интегралов при $\epsilon \to 0$
изучается, напр., в [4],гл.I.

Сравним скорости сходимости для индивидуальных функций при разных $\alpha>0$ в
зависимости от $\epsilon$.

При любых $\beta>\alpha$ для $f\in L_p,p\ge 1$ при любом $\epsilon>0$
$$
\|f(\cdot)-\frac {1}{(2\pi)^d}\int_{\mathbb R^d}\widehat{\phi_\beta} (y)f(\cdot-\epsilon
y)dy\|_p\le c(\alpha,\beta)\|f(\cdot)-\frac {1}{(2\pi)^d}\int_{\mathbb
R^d}\widehat{\phi_\alpha} (y)f(\cdot-\epsilon y)dy\|_p.
$$
Для доказательства применяем лемму 1.

Функция
$$
\psi(x)=\frac {1-\phi_\beta (x)}{1-\phi_\alpha (x)}\in W_1,
$$
напр., в силу того, что $\psi-1\in C_0(\mathbb R^d)\cap C^\infty(\mathbb R^d \setminus
\{0\})$ и при $N>\frac d2$ все производные порядка не выше $N$ принадлежат $L_2(R^d)$
(см. [3],6.5).

\emph{II. Сравнение дифференциальных операторов с постоянными коэффициентами.}

При $d=1$ задача о существовании неравенства
$$
\|Q(-i\frac {d}{dx})f\|_q\le c\|P(-i\frac {d}{dx})f\|_p,
$$
где $Q$ и $P$ - многочлены, исследована полностью. Найдены три разных
критерия, т.е. необходимые и достаточные условия одновременно, для функций на окружности
$\mathbb T$, полуоси $\mathbb R_+$ и оси $\mathbb R$ [7].

Здесь рассмотрим вопрос о существовании неравенства вида
$$
\|Q(-i\frac {d}{dx})f\|_q\le c(\|P_1(-i\frac {d}{dx})f\|_{p_1} + \|P_2(-i\frac
{d}{dx})f\|_{p_2})
$$
на прямой.

Для этого сумму операторов $P_1$ и $P_2$ представим в виде суммы двух свёрток с мерами.

\underline{Лемма 2}.

Если $deg Q\le r=deg P_1$ и $deg P_2\le r$, а общие вещественные нули $P_1$ и $P_2$, если
таковые имеются, являются нулями $Q$, то существуют функции $h_1$ и $h_2$ из $W_1(\mathbb
R)$ такие, что при всех вещественных $x$
$$
Q(x)=h_1(x)P_1(x)+h_2(x)P_2(x).
$$
Точнее, $h_1$ и $h_2$ удовлетворяют условию  $Lip 1$ на оси, и  $h_1(x)=\frac {Q(x}{P_1(x)}$ 
при больших $|x|$, а $h_2(x)=0$.

Доказательство.

Если $P_1(x)\ne 0$ на оси, то полагаем $h_1(x)=\frac {Q(x)}{P_1(x)}$ и $h_2(x)=0$. В
противном случае вне некоторых достаточно малых окрестностей вещественных нулей $P_1$, в
которых отсутствуют нули $P_2$, не являющиеся нулями $P_1$, полагаем $h_1(x)=\frac
{Q(x)}{P_1(x)}$ и $h_2(x)=0$. А в каждой из этих окрестностей полагаем $h_2(x)=\frac
{Q(x)-h_1(x)P_1(x)}{P_2(x)}$, где $h_1$ линейная, напр., определяемая значениями в концах
окрестности.

Лемма доказана.

Если $f$ и  $f^{(r)} \in L_1(\mathbb R)$, напр. (предполагается, конечно, что $f^{(r-1)}$
локально абсолютно непрерывна), то $\lim f^{(\nu)}(x)=0$ при $|x|\to \infty$ и $0\le
\nu\le r-1$. В этом случае
$$
\widehat{P_1(-i\frac {d}{dx})f} =P_1(y)
$$
и, значит,
$$
\widehat {Q(-i\frac {d}{dx})f}= Q(y)=h_1(y)\widehat{ P_1(-i\frac {d}{dx})f}
+h_2(y)\widehat {P_2(-i\frac {d}{dx})f}.
$$
Но $h_1=\widehat {d\mu_1}$ и $h_2=\widehat{ d\mu_2}$. Поэтому почти всюду
$$
Q(-i\frac {d}{dx})f=P_1(-i\frac {d}{dx})f*d\mu_1+P_2(-i\frac {d}{dx})f*d\mu_2,
$$
откуда и следует искомое неравенство. Если $degQ =r$, то $p_1=q$, а поскольку $\widehat
{{h_2}}\in L_s$ при любом $s\ge1$, то в силу неравенства Юнга $p_2$ - любое число из
$[1,q]$. Если же $degQ<r$, то по той же причине ($\widehat{h_1}\in L_s,s\ge 1$) и
$p_1$ - любое из $[1,q]$.

В случае $d\ge 2$ алгебра многочленов другая, но и в этом случае можно продвинуться
подобным образом и не только для операторов, близких к эллиптическим.

\bigskip

Литература

\bigskip

[1] Б. М. Макаров, А. Н. Подкорытов, Лекции по вещественному анализу. Санкт-Петербург,
БХВ-Петербург, 2011.

[2] Р. М. Тригуб, \emph{Абсолютная сходимость интегралов Фурье, суммируемость рядов
Фурье и приближение полиномами функций на торе}, Изв. АН СССР, с.м., 44:6 (1980),
1378--1408.

[3] E. Liflyand, S. Samko and R. Trigub,  \emph{The Wiener algebra of absolutely
convergent Fourier integrals: an overwiew}. Analysis and Math. Physics, Springer, 2:1
(2012), 1--68.

[4] E. M. Stein and G. Weiss, Introduction to Fourier Analysis on Euclidean
Spaces. Princeton Univ.Press, Princeton, 1971.

[5] H. S. Shapiro, \emph{Some Tauberian theorems with applications to approximation
theory}, Bull. Amer. Math. Soc. 74:3 (1968), 500--504.

[6] R. Trigub, E. Belinsky, \emph{Fourier Analysis and Approximation of Functions},
Kluwer-Springer, 2004.

[7] Р. М. Тригуб, \emph{О сравнении линейных дифференциальных операторов},
Матем.заметки, 82:3 (2007), 426--440

\bigskip

Донецкий национальный университет (Украина)

 roald.trigub@gmail.com

\end{document}